\documentclass[14pt]{extarticle}
\usepackage{extsizes}
\usepackage{amssymb,amsfonts,amsthm}
\usepackage{amsmath}
\usepackage[makeroom]{cancel}
\usepackage{mathtools}
\usepackage{mathrsfs,pifont}
\usepackage{slashed,mathabx} 
\usepackage[bbgreekl]{mathbbol}

\usepackage[top=1in, bottom=1.25in, left=0.75in,right=0.75in]{geometry}

\usepackage[font={small,sl}]{caption}
\usepackage[all]{xy}

\usepackage{hyperref}
\usepackage[backrefs,msc-links]{amsrefs}

\newcommand{\C}{\mathbb{C}}

\newcommand{\Z}{\mathbb{Z}}
\newcommand{\Q}{\mathbb{Q}}
\newcommand{\R}{\mathbb{R}}
\newcommand{\N}{\mathbb{N}}

\newcommand{\bbG}{\mathbb{G}}

\newcommand{\bB}{\mathsf{B}}

\newcommand{\bF}{\mathsf{F}}

\newcommand{\bT}{\mathsf{T}}
\newcommand{\bZ}{\mathsf{Z}}
\newcommand{\bZb}{\overline{\mathsf{Z}}}
\newcommand{\bAb}{\overline{\mathsf{A}}}

\newcommand{\bA}{\mathsf{A}}

\newcommand{\bP}{\mathbb{P}}

\newcommand{\cK}{\mathscr{K}}
\newcommand{\cL}{\mathscr{L}}

\newcommand{\cO}{\mathscr{O}}

\newcommand{\cU}{\mathscr{U}}

\newcommand{\cM}{\mathscr{M}}
\newcommand{\Mbar}{\overline{\cM}}

\newcommand{\cC}{\mathscr{C}}
\newcommand{\sS}{\scr{S}}

\newcommand{\yY}{\mathscr{Y}\!}

\newcommand{\cI}{\mathscr{I}}

\newcommand{\bal}{\boldsymbol{\alpha}}

\newcommand{\fg}{\mathfrak{g}} 
\newcommand{\fgh}{\widehat{\mathfrak{g}}} 
 
\newcommand{\fh}{\mathfrak{h}}

\newcommand{\Taut}{\mathsf{Taut}}

\newcommand{\rdd}{/\!\!/}
\newcommand{\Ct}{\mathbb{C}^\times}

\newcommand{\fC}{\mathfrak{C}}

\newcommand{\Hd}{{H}^{\raisebox{0.5mm}{$\scriptscriptstyle \bullet$}}}
\newcommand{\Hhd}{{H}_{\raisebox{0.5mm}{$\scriptscriptstyle \bullet$}}}
\newcommand{\bSd}{{\mathsf{S}}^{\raisebox{0.5mm}{$\scriptscriptstyle
      \bullet$}}}

\newcommand{\glh}{\widehat{\mathfrak{gl}}} 
\newcommand{\gl}{\mathfrak{gl}} 
 
\newcommand{\glhh}{\widehat{\widehat{\mathfrak{gl}}}}

\DeclareMathOperator{\Def}{Def}

\DeclareMathOperator{\Ob}{Ob}
\DeclareMathOperator{\Lie}{Lie}

\DeclareMathOperator{\Hilb}{Hilb}

\DeclareMathOperator{\End}{End}

\DeclareMathOperator{\Hom}{Hom}
\DeclareMathOperator{\Aut}{Aut}

\DeclareMathOperator{\ev}{ev}
\DeclareMathOperator{\Pic}{Pic}

\DeclareMathOperator{\Stab}{Stab}
\DeclareMathOperator{\Attr}{Attr}
\DeclareMathOperator{\diag}{diag}
\DeclareMathOperator{\codim}{codim}

\DeclareMathOperator{\Db}{D^b}

\newcommand{\vir}{\textup{vir}}

\newcommand{\pt}{\textup{pt}}

\newcommand{\cF}{\mathscr{F}}

\newcommand{\tO}{\widehat{\mathscr{O}}}
\newcommand{\hX}{\widehat{X}}

\DeclareMathOperator{\Coker}{Coker}
\DeclareMathOperator{\tr}{tr}

\DeclareMathOperator{\Coh}{Coh}

\newcommand{\bo}[1]{\textbf{(}\textsf{#1}\textbf{)}}

\newtheorem{Theorem}{Theorem}

\theoremstyle{definition}

\makeatletter
\def\blfootnote{\xdef\@thefnmark{}\@footnotetext}
\makeatother

\begin{document}

\title{On the crossroads of enumerative geometry and geometric representation theory} 
\author{Andrei Okounkov}
\date{} 
\maketitle



\noindent
The subjects in the title are interwoven in many different and very
deep ways. I recently wrote several expository accounts
\cites{PCMI,slc,takagi} 
that reflect a certain range of developments, but even in
their totality they cannot be taken as a comprehensive survey. 
In the format of a $30$-page contribution aimed at a general mathematical
audience, I have decided to illustrate some of the basic ideas in one
very interesting example --- that of $\Hilb(\C^2,n)$, hoping to spark 
the curiosity of colleagues in those numerous fields of study where
one should expect applications.

\blfootnote{It is impossible to overestimate the 
influence of my friends and collaborators M.~Aganagic,
  R.~Bezrukavnikov, P.~Etingof, I.~Loseu, D.~Maulik, N.~Nekrasov,
  R.~Pandharipande on the thoughts and ideas presented here. 
I am very grateful to the Simons Foundation Simons Investigator program,  the Russian Academic
Excellence Project ’5-100’, and the RSF grant 16-11-10160 for their 
support of my research.}

\section{The Hilbert scheme of points in $\C^2$} 

\subsection{Classical geometry}

\subsubsection{} An unordered collection $P=\{p_1,\dots,p_n\}$ of
distinct points in the plane $\C^2$ is uniquely specified by the
corresponding ideal 
$$
I_P = \{f(p_1) = \dots = f(p_n) =0\} \subset \C[x_1,x_2]
$$
in the coordinate ring $\C^2$. The codimension of this ideal, 
i.e. the dimension of the quotient 
$$
\C[x_1,x_2]/ I_P = \textup{functions on $P$}
\,\overset{\textup{def}}=\, \cO_P \,,
$$
clearly equals $n$. 

If the points $\{p_i\}$ merge, the
limit of $I_P$ stores more information than just the location of
points. For instance, for two points, it remembers the direction along 
which they came together. One defines 
$$
\Hilb(\C^2,n) = \{ \textup{ideals  }I \subset \C[x_1,x_2] \textup{
  of codimension $n$}\} \,,
$$
and with its natural scheme structure \cite{FDA,Koll} this turns out to be a smooth 
irreducible algebraic variety --- a special feature of the Hilbert
schemes of surfaces that fails very badly in higher dimensions. 

 The map $I_P \mapsto P$ extends to a natural map 
 \begin{equation}
\pi_{\Hilb}: \Hilb(\C^2,n) \to (\C^2)^n / S(n) \label{piHilb}
\end{equation}
which is proper and birational, in other words,
a resolution of singularities of $(\C^2)^n / S(n)$. This makes
$\Hilb(\C^2,n)$ an instance of an \emph{equivariant symplectic 
resolution} ---  a very special class of algebraic varieties
\cite{Kal,Beau} that
plays a central role in the current development of both 
enumerative geometry and geometric representation theory.
This general notion axiomatizes two key features of $\pi_{\Hilb}$: 
\begin{itemize}
\item the source of $\pi$ is an algebraic symplectic
  variety (here, with the symplectic form induced from that of
  $\C^2$), 
\item the map is equivariant for an action of a torus $\bT$ 
that contracts the target to a point (here, $\bT$ are the diagonal 
matrices in $GL(2)$ and the special point is the origin in
$(\C^2)^n$). 
\end{itemize}

\subsubsection{}

Both enumerative geometry and geometric representation theory really 
work with algebraic varieties $X$ and correspondences, that is, 
cycles (or sheaves, etc.) in $X_1 \times X_2$, considered up to a
certain equivalence. These one can compose geometrically and they 
form a nonlinear analog of matrices of linear algebra and classical 
representation theory. To get to vector spaces and matrices, one 
considers functors like
 the equivariant cohomology $\Hd_\bT(X)$, the equivariant 
K-theory etc., with the induced action of the correspondences. 

Working with equivariant cohomology whenever there is a torus action 
available is highly recommended, in particular, because: 
\begin{itemize}
\item equivariant cohomology is in many ways simpler than
  the 
  ordinary, while also more general. E.g.\ the spectrum of the ring 
$\Hd_\bT(\Hilb(\C^2,n))$ is a union of explicit essentially linear
subvarieties over all partitions of $n$. 
\item the base ring $\Hd_\bT(\pt,\Q)= \Q[\Lie \bT]$ of
  equivariant cohomology introduces parameters in the theory, on which
  everything depends in a very rich and informative way, 
\item equivariance is a way to trade global geometry for
  local parameters. For instance, all formulas in the classical (that
  is, not quantum) geometry of Hilbert schemes of points generalize 
\cite{EGL} to
  the general surface $\sS$ with the substitution
$$
c_1(\sS) = t_1 + t_2 \,, \quad c_2(\sS) = t_1 t_2  \,, 
$$
where $\diag(t_1,t_2) \in \Lie \bT$. 
\end{itemize}

\subsubsection{}
The fundamental correspondence between Hilbert schemes of points 
is the Nakajima correspondence 
$$
\alpha_{-k}(\gamma) \subset \Hilb(\sS,n+k) \times \Hilb(\sS,n) 
$$
formed for $k>0$ by pairs of ideals $(I_1,I_2)$ such that $I_1 \subset I_2$ and 
the quotient $I_2/I_1$ is supported at a single point located along a 
cycle $\gamma$ in the surface $\sS$. This cycle is Lagrangian if $\gamma
\subset \sS$ is an algebraic curve. For $\sS=\C^2$, there are 
very few $\bT$-equivariant choices for $\gamma$ and they are all 
proportional, e.g. 
$$
[\{x_2 = 0 \} ] = t_2 [\C^2], \quad [0] = t_1 t_2 [\C^2] \,. 
$$
One defines $\alpha_{k}(\gamma)$ for $k>0$ as $(-1)^k$ times the
transposed correspondence (see section 3.1.3 in \cite{MO} on the author's preferred
way to deal with signs in the subject). 

The fundamental result of Nakajima \cite{NakHart} is the following commutation
relation\footnote{
in which the commutator is the supercommutator for odd-dimensional 
cycles $\gamma,\gamma'$. Similarly, the symmetric algebra in 
\eqref{fock} is taken in the $\Z/2$-graded sense.}
\begin{equation}
  \label{Nakcomm}
  [ \alpha_{n}(\gamma), \alpha_m(\gamma') ] = 
- (\gamma, \gamma')  \, n \, \delta_{n+m} \,, 
\end{equation}
where the intersection pairing $ (\gamma, \gamma')$ 
for $\sS=\C^2$ is $([0],[\C^2]) =1$. 
One recognizes in \eqref{Nakcomm} the commutation relation for the 
Heisenberg Lie algebra $\glh(1)$ --- a central extension of the 
commutative Lie algebra of Laurent polynomials with values in 
$\gl(1)$. The representation theory of this Lie algebra is very
simple, yet very constraining, and one deduces the identification
\begin{equation}
\bigoplus_{n\ge 0} \Hd_\bT(\Hilb(\sS,n)) \cong \bSd 
\left( \textup{span of }\{\alpha_{-k}(\gamma)\}_{k>0}\right) 
\label{fock}
  \end{equation}
with the Fock module generated by the vacuum 
$\Hd_\bT(\Hilb(\sS,0))= \Hd_\bT(\pt)$. 

Fock spaces (equivalently, symmetric functions) are everywhere in
mathematics and mathematical physics and many 
remarkable computations and phenomena are naturally expressed in this
language. My firm belief is that geometric construction, in particular
the DT theory of 3-folds to be discussed below, are the best known way
to think about them. 

\subsubsection{}

The identification \eqref{fock} is a good example to illustrate the general idea that the best
way to understand an algebraic variety $X$ and, in particular, its 
equivariant cohomology $\Hd_\bT(X)$, is to construct interesting 
correspondences acting on it. 

For a general symplectic resolution 
\begin{equation}
  \label{sympl_res}
  \pi: X \to X_0 
\end{equation}
the irreducible components of the \emph{Steinberg} variety $X
\times_{X_0} X$ give important correspondences. For 
$X = \Hilb(\C^2,n)$, these will be quadratic in Nakajima
correspondences, and hence not as fundamental. In Section \ref{s_gaqg} we
will see one general mechanism into which $\alpha_{k}(\gamma)$ fit. 

\subsubsection{}\label{s_Lehn} 

Since $\glh(1)$ acts irreducibly in \eqref{fock}, it is natural to
express all other geometrically defined operators in terms of
$\alpha_k(\gamma)$. Of special importance in what follows will be 
the operator of cup product (and also of the quantum product) by 
Chern classes of the tautological bundle $\Taut= \C[x_1,x_2]/I$ over
the Hilbert scheme. 

The operator of multiplication by the divisor $c_1(\Taut)$ was computed 
by M.~Lehn in \cite{Lehn} as follows 
\begin{align}
  \label{Lehn}
  c_1(\Taut) \cup = &- \tfrac12 \sum_{n,m>0} 
\big(\alpha_{-n} \, \alpha_{-m}
  \, 
\alpha_{n+m} +\alpha_{-n-m}  \, \alpha_{n} \, \alpha_{m} \big) \notag \\
&+ 
(t_1 + t_2) \sum_{n>0} \frac{n-1}{2} \,\alpha_{-n}  \, \alpha_{n} 
\end{align}
with the following convention\footnote{Note that it differs by a sign
  from the convention used in \cite{MO}.} about the arguments of the $\alpha$'s. 
If, say, we have 3 alphas, and hence need 3 arguments, we take the K\"unneth 
decomposition of 
$$
[ \textup{small diagonal} ] = [\C^2] \times [0] \times [0]  
\in \Hd_\bT ( (\C^2)^3) \,.
$$
Similarly, $\alpha_{-n}  \, \alpha_{n}$ is short for $\alpha_{-n}([\C^2])
\, 
\alpha_{n}([0])$. 
With this convention, \eqref{Lehn} is clearly an operator of
cohomological degree $2$. There is a systematic way to prove formulas
like \eqref{Lehn} in the framework of Section \ref{s_gaqg}, see
\cite{SmirnovR}. 

\subsubsection{}
A remarkable observation, made independently by several people, is
that the operator \eqref{Lehn} is identical to the second-quantized 
trigonometric 
Calogero-Sutherland operator --- an classical object in many-body
systems and symmetric functions (its eigenfunctions being the Jack
symmetric polynomials), see e.g.\ \cite{CostGroj} for a comprehensive 
discussion. 

The quantum CS Hamiltonian\footnote{Note that the 
well-known, but still remarkable strong/weak duality 
$\theta \mapsto 1/\theta$ in the CS model becomes simply 
the permutation of the coordinates in $\C^2$ in the geometric 
interpretation.}
\begin{equation}
H_{CS} = \tfrac12 \sum_{i=1}^N \left(w_i \frac{\partial}{\partial
    w_i}\right)^2  + \theta(\theta-1) \sum_{i<j\le N}
\frac{1}{|w_i-w_j|^2} 
\,, \quad\theta = -\left( t_1 \big/ t_2\right)^{\pm 1}\,, \label{HCS}
\end{equation}
describes a system of $N$ identical particles interacting with
$|w|^{-2}$-potential on the unit circle $|w|=1$. After conjugation by 
an eigenfunction $\prod_{i<j} (w_i-w_j)^{\theta}$, it preserves
symmetric Laurent polynomials in $w_i$ and stabilizes as $N\to \infty$ to a
limit in which the left-movers and right-movers (that is, symmetric polynomials
in $w_i$ and polynomials in $w_i^{-1}$) decouple. This limit is
\eqref{Lehn} with $\alpha_{-k}$ proportional to multiplication by
$\sum w_i^k$. 

Another interpretation of the same equation \eqref{Lehn} 
is an integrable quantum version of the Benjamin-Ono equation
of 1-dimensional
hydrodynamics, see in particular \cite{AW}. The BO equation 
describes waves on a 1-dimensional surface of a fluid of 
infinite depth, and it involves the Hilbert transform --- a nonlocal
operation. This nonlocal operation is precisely responsible for 
the term $\sum_{n>0} n \, \alpha_{-n} \, \alpha_{n}$ which is present
in \eqref{Lehn} and looks a bit unconventional when expressed in terms of the 
field $\bal(\zeta) = \sum \alpha_{-n} \zeta^n$. Note that other terms in 
\eqref{Lehn} are the normally ordered 
constant terms in $\bal(\zeta)^3$ and $\bal(\zeta)^2$, 
respectively. 

These by now classical connections are only a preview of the kind of
connections that exists between enumerative problems and quantum 
integrable systems in the full unfolding of the theory. Crucial
insights into this connection were made in the pioneering work of 
Nekrasov and Shatashvili \cites{NS1,NS2}. 

\subsection{Counting curves in $\Hilb(\C^2,n)$}

\subsubsection{}

Enumerative geometry of curves in an algebraic variety $Y$ is a very
old subject in mathematics, with the counts like the 27 lines on a smooth 
cubic surface going as far back as the work of Cayley from 1849. While 
superficially the subject may be likened to counting points
of $Y$ over some field $\Bbbk$, the actual framework that the
geometers have to construct to do the counts looks very different from 
the number-theoretic constructions. In particular, the counts are
defined treating curves in $Y$ as an excess intersection
problem, with the result that the counts are invariant under
deformation even though there may be no way to deform actual curves. 

Also, the subject draws a lot of inspiration
from mathematical physics, where various curve counts are interpreted
as counts (more precisely, indices) of supersymmetric states in certain gauge or string
theories. This leaves a very visible imprint on the field, ranging from how
one organizes the enumerative data to what is viewed as an important 
goal/result in the subject. In particular, any given count, unless it
is something as beautiful as 27 lines, is viewed as only an
intermediate step in the quest to uncover universal structures that
govern a certain totality of the counts. 

\subsubsection{}

Depending on what is meant by a ``curve'' in $Y$, enumerative theories
come in several distinct flavors, with sometimes highly nontrivial
interrelations between them. My personal favorite among them is 
the Donaldson-Thomas theory of 3-folds \cite{DonTh,ThCass}, see, in particular, 
\cite{takagi} for a recent 
set of lecture notes. 

The DT theory views a curve $C\subset
Y$ as something defined by equations in $Y$, that is, as a subsheaf 
$\cI_C\subset \cO_Y$ of regular functions on $Y$ formed by those 
functions that vanish on $C$. The DT moduli space for $Y$ is thus 
\begin{equation}
\Hilb(Y,\textup{curves}) = \bigsqcup_{\beta,n}
\Hilb(Y,\beta,n)\label{Hilb-curves}
\end{equation}
where $\beta\in H_2(Y,\Z)_\textup{eff}$ is the degree of the curve and 
$n=\chi(\cO_C)$ is the holomorphic Euler characteristic of 
$\cO_C= \cO_Y/\cI_C$. In particular, if $C$ is a smooth connected curve of 
genus $g$ then $n=1-g$. 

\begin{figure}[!htbp]
  \centering
   \includegraphics[scale=0.59]{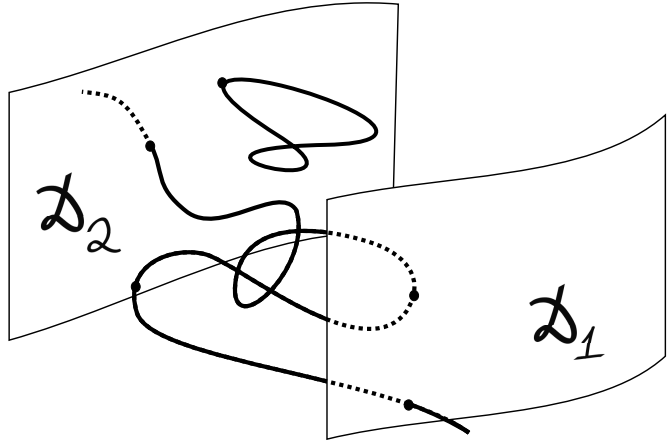}
 \caption{A fundamental object to count in DT theory are 
algebraic curves, or more precisely, subschemes $C\subset Y$ of given
$(\beta,\chi)$ constrained by how they meet a fixed divisor $D\subset Y$.}
\label{f_DT}
\end{figure}

While similar in construction and
universal properties to the Hilbert schemes of point in surfaces, 
Hilbert schemes of 3-folds are, in general, highly singular 
varieties and nearly nothing is known about their dimension or 
irreducible components. However, viewed as moduli of sheaves of on $Y$ of 
the form $\{\cI_C\}$, they have a good deformation theory and 
thus a virtual fundamental cycle of (complex) dimension $\beta \cdot c_1(Y)$ \cite{ThCass}. The DT curve counts are
defined by pairing this virtual cycle against natural cohomology
classes, such as those pulled back from $\Hilb(D,\textup{points})$ via
a map that assigns to $C$ its pattern of tangency to a fixed smooth divisor 
$D\subset Y$, see Figure 
\ref{f_DT} and the discussion in Section 2 of \cite{takagi}. 
%
Note that the divisor $D=\bigsqcup D_i$ may be disconnected,
as in Figure \ref{f_DT}, 
in which case curve counts really define a tensor in a tensor product of 
several Fock spaces, e.g. an operator from one Fock space to another if
there are two components. I believe that this is 
the  geometric source of great many, if not all,
interesting tensors in Fock spaces. 

There is a broader world of DT counts, in which one counts other
1-dimensional sheaves on 3-folds (a very important example being the 
Pandharipande-Thomas moduli spaces of stable pairs \cite{PT1}), and even
more broadly stable objects in other categories that look like
coherent sheaves on a smooth 3-fold. For all of these, the deformation
theory has certain self-duality features that makes K-theoretic and
otherwise refined enumerative information a well-behaved and a very
interesting object to study. 

\subsubsection{}\label{szn} 

Since the virtual dimension does not depend on $n$, it is convenient
to organize the DT counts by summing over all $n$ with a weight $z^n$,
where $z$ is a new variable. These are conjectured to be rational
functions of $z$ with poles at roots of unity \cites{MNOP1,MNOP2}. This is known
in many important cases \cites{MOOP,PP2,SmirnovRat,Toda1} and may be put into a larger 
conjectural framework as in \cite{MDT}, see Section \ref{s_Kcounts}.


\subsubsection{}\label{sTQFT} 

\begin{figure}[!htbp]
  \centering
   \includegraphics[scale=0.62]{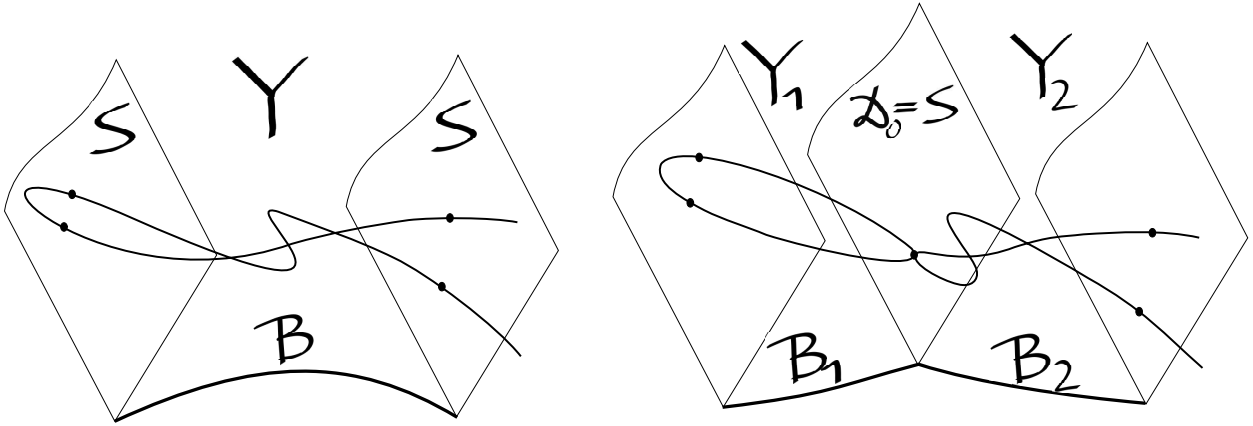}
 \caption{Basic building blocks of DT theory are $\sS$ bundles over
   a curve $B$ as in \eqref{fibrS}. As $B$ degenerates to a nodal
   curve $B_1 \cup B_2$, $Y$ degenerates to a transverse union of $Y_1$ and $Y_2$
   along a smooth divisor $D_0\cong\sS$. DT counts satisfy a gluing formula
   for such degenerations.}
\label{f_deg}
\end{figure}

An algebraic analog of cutting $Y$ into pieces is a degeneration of 
$Y$ to a transverse union of $Y_1$ and $Y_2$ along a smooth divisor
$D_0$ as in Figure \ref{f_deg}. A powerful result of 
Levine and Pandharipande \cite{LP} shows that any smooth projective
3-fold can be linked to a product of projective spaces (or any other 
basis in algebraic cobordism) by a sequence of such moves. 
The DT 
counts satisfy a certain gluing formula for such degeneration 
\cite{LiWu} in which the divisor $D_0$ is added to the 
divisors from Figure \ref{f_DT}.  This highlights the importance of understanding the DT counts 
in certain basic geometries which can serve as building blocks for 
arbitrary 3-folds. 

One set of basic geometries is formed by  $\sS$-bundles over a curve $B$
\begin{equation}
  \label{fibrS}
\xymatrix{
  \sS \ar@{^{(}->}[r]  & Y
\ar[d]_g \\
& B \,,
} 
\end{equation}
where $\sS$ is
a smooth surface\footnote{It suffices to take $\sS\in
  \{A_0=\C^2,A_1,A_2\}$, where $A_n$ is the minimal resolution of the
  corresponding surface singularity,  to generate a basic set of
  counts.}. 
One take $D=\bigcup g^{-1}(b_i)$ for $\{b_i\} \subset B$ and, by 
degeneration, this defines a TQFT on $B$ with the space of states
$\Hd_\bT(\Hilb(\sS))$, where 
$\bT\subset \Aut(\sS)$ is a maximal torus. This TQFT structure
 is captured by the counts for $B=\bP^1
\supset \{b_1,b_2,b_3\}$, which define a new, $z$- and 
$\beta$-dependent 
supercommutative multiplication in this Fock space. It is a
very interesting question to describe this multiplication 
explicitly\footnote{In particular, there are very interesting results
  and conjectures for K3 surface fibrations, see \cite{Ober}. Note that 
the dimension counts work out best when $c_1(\sS)=0$.}. 

Geometric representation theory provides an answer when $\sS$ itself is 
a symplectic resolution, which concretely means an ADE surface --- a 
minimal resolution of the corresponding surface singularity $\sS_0$. In this
case, $\Hilb(\sS,n)$ is a symplectic resolution of $(\sS_0)^n/S(n)$,
just like \eqref{piHilb}, and the multiplication above is
nothing but the quantum multiplication in $\Hd_\bT(\Hilb(\sS))$, see
\cites{MauObl,OP4}. 

\subsubsection{}
Recall that the quantum product $\star$ is a supercommutative, 
associative deformation of the classical $\cup$-product in
$\Hd_\bT(X)$ whose structure constants are the counts of 3-pointed 
rational curves in $X$ 
\begin{equation}
  \label{star}
  (\alpha \star \beta, \gamma) = 
\sum_{d\in H_2(X,\Z)_{\textup{eff}}} z^d \quad 
\begin{matrix}
\textup{$\#$ of degree $d$ rational curves}\\
\textup{meeting cycles dual to $\alpha, \beta, \gamma$} \,. 
\end{matrix}
\end{equation}
This is made mathematically precise using the correspondence 
\begin{equation}
\sum_d z^d \ev_* \left(\left[\Mbar_{0,3}(X,d)\right]_\textup{vir}\right) 
\in \Hhd(X \times X \times X)[[z]] \label{cM03}
\end{equation}
obtained from the virtual fundamental cycle of the moduli space of 
3-pointed stable rational maps to $X$, see \cite{mirror_book} for an
introduction. 

A closely related structure is the \emph{quantum differential
  equation}, or Dubrovin connection $\nabla_X$, which is a flat connection on trivial bundle over 
$H^2(X,\C)\owns \lambda$ with fiber $\Hd_\bT(X,\C)$. Its flat sections satisfy 
\begin{equation}
\frac{d}{d\lambda} \Psi(z) = \lambda \star \Psi(z) \,, \quad 
\frac{d}{d\lambda}  z^d = (\lambda,d) \, z^d \,, \label{qde}
\end{equation}
and contain very important enumerative information. For 
$X=\Hilb(\C^2)$, the quantum multiplication ring is generated by
the divisor, so the two structures are really the same. 

\subsubsection{}\label{s_pq}

With the $z\mapsto -z$ substitution, the $\star$-deformation of
\eqref{Lehn} for $X=\Hilb(\C^2)$ was computed in 
\cite{OP4} as follows 
\begin{equation}
  \label{starLehn}
  c_1(\Taut) \star = c_1(\Taut) \cup + (t_1+t_2) 
\sum_{d>0} d \frac{z^d}{1-z^d} \, \alpha_{-d} \, \alpha_d + \dots\,, 
\end{equation}
with the same convention about the arguments of $\alpha_n$ as in 
\eqref{Lehn} and with
dots denoting a scalar operator of no importance for us now. 
Note the simplicity of the purely quantum terms\footnote{
If one interprets \eqref{Lehn} as a quantum version of the Benjamin-Ono 
equation then the new terms deform it to a quantization of the 
intermediate long wave (ILW) equation. This observation has been 
rediscovered by many authors.}. 

The largest and the richest class of equivariant symplectic 
resolutions known to date is formed by the Nakajima quiver 
varieties \cite{Nakq1}, of which $\Hilb(\C^2)$ is an example. 
Formula \eqref{starLehn} illustrates many general features of
the quantum cohomology of Nakajima quiver varieties proven in 
\cite{MO}, such as: 
\begin{itemize}
\item  the purely quantum terms are given by 
a rational function with values in Steinberg correspondences
  \begin{equation}
    \label{pq1}
    \textup{purely quantum} \in \hbar\,  H_{\textup{top}}(X \times_{X_0}
    X)
\otimes \Q(z) \,, 
  \end{equation}
where 
$$
\hbar = - (t_1+t_2)
$$
 is the equivariant weight of the symplectic
form. 
\item the shift $z\mapsto -z$ is an example of the shift by 
a canonical element of $H^2(X,\Z/2)$, called the 
theta-characteristic in \cite{MO}. 
\item there is a certain Lie algebra $\fg_Q$ associated to 
an arbitrary quiver $Q$ in \cite{MO}, whose positive roots are represented
by effective curve classes 
$d\in H_2(X,\Z)$. Among these, 
there is a finite set of \emph{K\"ahler roots} of $X$ such that 
  \begin{equation}
    \label{pq2}
    \lambda\star = \lambda \cup 
- \, \hbar \!\sum_{d\in \{\textup{positive roots}\}} (\lambda,d) \, \frac{z^d}{1-z^d}  \, C_d 
\,, 
  \end{equation}
where $C_d \in \fg_{-d} \, \fg_d$ is the corresponding root component
of the Casimir element, that is, the image of the invariant 
bilinear form on $\fg_{-d} \otimes \fg_{d}$. 
For $X=\Hilb(\C^2,n)$, the quiver is the quiver with one vertex and
one loop with $\fg_Q = \glh(1)$, $C_d = \alpha_{-d} \, \alpha_d$, and 
\begin{equation}
\textup{positive K\"ahler roots}  = \{1,\dots,n\} \subset \Z
\,. \label{KR_hilb}
\end{equation}
\item the Cartan subalgebra $\fh \subset \fg_Q$ acts by 
central elements and by the 
ranks of the tautological bundles.  The operator of cup products by 
other characteristic classes of the tautological bundles, together
with $\fg$, generate a Hopf algebra deformation $\yY(\fg_Q)$ of $\cU(\fg_Q[t])$
known as the Yangian. 

The operators of quantum multiplication form a remarkable family 
of maximal commutative subalgebras of $\yY(\fg_Q)$
known as the Baxter subalgebras in the theory of quantum integrable 
systems, see Section \ref{sBaxt}. They are parametrized by $z$ and,  as $z\to 0$, they become 
the algebra 
$\yY(\fh) \subset \yY(\fg_Q)$ of cup products by tautological classes. 

The identification between the $\star$-product ring and Baxter's
quantum integrals of motion was predicted by Nekrasov and Shatashvili
based on their computation of the spectra of the operators. This served as very important inspiration for \cite{MO}. 
\item
The Yangian description identifies the quantum differential equation 
with the \emph{Casimir connection} for the Lie algebra $\fg_Q$, as
studied (in the finite-dimensional case) in \cite{TolCas}. 
This fits very
nicely 
with the conjecture of Bezrukavnikov and collaborators about the 
monodromy $\nabla_X$, see below, and was another important 
inspiration for \cite{MO}, see the historical notes there. 
\end{itemize}

For general symplectic resolutions, there is a definite gap between what 
is known abstractly, and what can be seen in known
examples. I expect that a complete classification of the equivariant 
symplectic resolutions is within the reach of the current generation
of algebraic geometers, and we will see how representative the known 
examples are.  For general symplectic
resolutions,  the Steinberg correspondence in 
\eqref{pq1} are constructed in \cite{BMO}, while the rationality in
$z$ remains abstractly a conjecture that can be checked in all known
cases. 


\subsubsection{}

A generalization of these structures appears 
in enumerative K-theory. There, instead of pairing virtual cycles with cohomology
classes, we compute the Euler characteristics of
natural sheaves, including the virtual structure sheaf
$\cO_\vir$, on the moduli spaces in question (here, the Hilbert
scheme 
of curves in $Y$). 

{} From its very beginning, K-theory has been inseparable from the 
indices of differential operators and related questions in
mathematical physics. Equivariant 
K-theoretic DT counts represent a Hamiltonian 
approach to supersymmetric indices in a certain physical theory (namely, the theory on a D6
brane), in which the space is $Y$ and the time is 
periodic\footnote{Or, more precisely, quasiperiodic, with a twist  
by an element of $\bT\subset \Aut(Y)$ after a full circle of
time.}. 
Morally, what 
one computes is the index of a certain infinite-dimensional Dirac
operator as a representation of $\Aut(Y)$ which is additionally 
graded by $(\beta,n)$. Because this is the index of a Dirac operator, 
the right analog of the virtual cycle is the symmetrized virtual 
structure sheaf
\begin{equation}
\tO_\vir = \cO_\vir \otimes \cK^{1/2}_\vir \otimes \dots \label{tOvir}
\end{equation}
where $\cK^{1/2}_\vir$ is a square root of the virtual canonical
bundle $\cK_\vir$, the importance of which was emphasized by Nekrasov
\cite{Zth}, and the existence of which is shown in \cite{MDT}. The dots in 
\eqref{tOvir} denote a certain further twist by a tautological line
bundle of lesser importance, see \cite{MDT}. 

\subsubsection{}

While it is not uncommon for different moduli spaces to give the same
or equivalent cohomological counts, the K-theoretic counts really feel
every point in the moduli space and are very sensitive to the exact 
enumerative setup. In particular, for both the existence of
\eqref{tOvir} and the computations with this sheaf, certain
self-duality features of the DT deformation theory are crucial. It
remains to be seen whether computations with moduli spaces like 
$\Mbar_{0,3}(X)$, that lack such self-duality, can really reproduce
the K-theoretic DT counts. 

In K-theory, the best setup for counting curves in $X=\Hilb(\C^2,n)$
is the moduli space of stable \emph{quasimaps} to $X$, see \cite{CKM}. 
Recall that
$$
X = \left\{ x_1,x_2 \in \End(\C^n), v \in \C^n \big| [x_1,x_2]=0
\right\}
\rdd GL(n)\,, 
$$
where the stability condition in the GIT quotient is equivalent to 
$\C[x_1,x_2] v$ spanning $\C^n$. By definition \cite{CKM}, a stable quasimap from $B$ to
a GIT quotient is a 
map to the quotient stack\footnote{i.e.\ a principal $G$-bundle on $B$
together with a section of the associated bundle of prequotients,
where $G$ is the group by which we quotient. For $G=GL(n)$, a
principal $G$-bundle is the same as a vector bundle on $B$ of rank
$n$.
} that
evaluates to a stable point away from a finite set of point in
$B$. In return for allowing such singularities, quasimaps offer
many technical advantages. 

If $Y$ in \eqref{fibrS} is a fibration in $\sS=\C^2$, one can consider 
quasimap sections of the corresponding $X$-bundle over $B$, and these 
are easily seen to be identical to the Pandharipande--Thomas stable
pairs for $Y$. Recall that by definition \cite{PT1}, a stable pair is a complex 
of the form 
$$
\cO_Y \xrightarrow{\quad s \quad} \cF
$$
where $\cF$ is a pure 1-dimensional sheaf and $\dim \Coker s = 0$. 
For our $Y$, $g_* \cF$ is a vector bundle on $B$, the section $s$
gives $v$, while $x_1,x_2$ come from multiplication by the coordinates
in the fiber. If the fiber $\sS$ contains curve, the picture
becomes modified and the PT spaces for the $A_n$-fibrations, $n>0$, are 
related to quasimaps via a certain sequence of wall crossings. 

\subsubsection{}
The K-theoretic quasimap counts to $\Hilb(\C^2)$ and, in fact, to all
Nakajima varieties have been computed in \cite{PCMI,OS}, and their
 structure  
is a certain $q$-difference deformation of what we have seen for the 
cohomological counts. In particular, the Yangian $\yY(\fg)$ is replaced
by a quantum loop algebra $\cU_\hbar(\fgh)$ formed by K-theoretic 
analogs of the correspondences that define the action of $\yY(\fg)$, 
see Sections \ref{s_Stab} and \ref{s_StabK}. 

\section{Geometric actions of quantum groups}
\label{s_gaqg}

Geometric representation theory, in the sense of making interesting
algebras act by correspondences, is a mature subject and its exposition
in \cite{CG} is a classic. Geometric construction of representation of 
quantum groups has been a very important stimulus in the development 
of the theory of Nakajima quiver varieties \cite{Nakq2, Nak3}. 

Below we discuss a 
complementary approach of \cite{MO}, which mixes geometry and algebra in
a different proportion. It has certain convenient hands-off features, in the sense that it constructs a certain category of 
representations without, for example, a complete description of the
algebra by generators and relations. In algebraic geometry, one certainly prefers 
having a handle on the category $\Coh(X)$ to a complete list of
equation that cut out $X$ inside some ambient variety, so the
construction should be of some appeal to algebraic geometers. It also 
interact very nicely with the enumerative question, as it has certain
basic compatibilities built in by design.

\subsection{Braiding} 

\subsubsection{} 
If $G$ is a group then the category of $G$-modules over a field 
$\Bbbk$ has a tensor
product --- the usual tensor product $M_1 \otimes M_2$ of vector
spaces in which an element $g\in G$ acts by $g\otimes g$. There is 
also a trivial representation $g\mapsto 1\in \End(\Bbbk)$, which is
the identity for $\otimes$. This
reflects the existence of a coproduct, that is, of an algebra homomorphism 
\begin{equation}
\Bbbk G \owns g  \xrightarrow{\quad \Delta \quad} 
g \otimes g \in \Bbbk G  \otimes \Bbbk G \,,\label{gg}
\end{equation}
where $\Bbbk G$ is the group algebra of $G$, with the counit 
$$
\Bbbk G \owns g \xrightarrow{\quad \varepsilon \quad} 1 \in \Bbbk \,.
$$
There are also dual module $M^*=\Hom_\Bbbk(M,\Bbbk)$ in
which $g$ acts by $(g^{-1})^T$, reflecting the antiautomorphism
$$
\Bbbk G \owns g  \xrightarrow{\quad \textup{antipode} \quad}  g^{-1}
\in \Bbbk G \,. 
$$
Just like the inverse in the group, the antipode is unique if it
exist, and so it will be outside of our focus in what follows, see
\cite{Ebook2}. 

An infinitesimal version of this for a Lie or algebraic group $G$ is to replace 
$\Bbbk G$ by the universal enveloping algebra $\cU(\fg)$,  
$\fg = \Lie G$, with the coproduct obtained from \eqref{gg} by 
Leibniz rule 
$$
\Delta(\xi) = \xi \otimes 1 + 1 \otimes \xi \,, \quad \xi \in \fg \,. 
$$
The multiplication, comultiplication, unit, counit, and the antipode
form an beautiful algebraic structure known as a Hopf algebra, see 
e.g.\ \cite{Ebook2}. Another classical example is the algebra $\Bbbk[G]$ of
regular functions on an algebraic group. Remarkably, the axioms of a Hopf
algebra are self-dual under taking duals and reversing all arrows. 
Observe that all of the above examples are either commutative, or
cocommutative. 

\subsubsection{}\label{s_Rmat}

Broadly, a quantum group is a deformation of the above examples in the
class of Hopf algebras. Our main interests is in Yangians and quantum 
loop algebras that are deformations of 
\begin{equation}
\cU(\fg[t]) \rightsquigarrow \yY(\fg) \,, \quad
\cU(\fg[t^{\pm 1}]) \rightsquigarrow \cU_\hbar(\fgh) \,,\label{defdef}
\end{equation}
respectively. Their main feature is the loss of cocommutativity. In
other words, the order of tensor factors now matters and 
$$
M_1 \otimes M_2 \not \cong M_2 \otimes M_1 \,, 
$$
in general, or at least the permutation of the tensor factors is no
longer an intertwining operator. While tensor categories are very familiar to algebraic geometers, this
may be an unfamiliar feature. But, as the representation-theorists
know, a mild noncommutativity of the tensor product makes the 
theory richer and more constrained. 

The Lie algebras $\fg[t]$ and $\fg[t^{\pm 1}]$ 
in \eqref{defdef} are $\fg$-valued functions on the 
additive, respectively multiplicative, group of the field and they
have natural automorphisms
$$
t \mapsto t+a \,, \quad \textup{resp.} \quad t \mapsto a t \,,
\quad a \in \bbG\,, 
$$
where we use $\bbG$ as a generic symbol for either an additive or 
multiplicative group. The action of $\bbG$ 
will deform to an automorphism of the quantum
 group\footnote{In fact, for $\cU_\hbar(\fgh)$, it is natural to view this
\emph{loop rotation} automorphism as part of the Cartan torus.} and we 
denote by $M(a)$ the module $M$, with the action precomposed by
an automorphism from $\bbG$. 

The main feature of the theory is the existence of intertwiner (known
as the braiding, or the $R$-matrix) 
\begin{equation}
  \label{Rchek}
  R^\vee(a_1-a_2): M_1(a_1) \otimes M_2(a_2) \to M_2(a_2) \otimes
  M_1(a_1) 
\end{equation}
which is invertible as a rational function of $a_1-a_2\in\bbG$ and develops a 
kernel and cokernel for those values of the parameters where the 
two tensor products are really not isomorphic. One often works with 
the operator $ R = (12) \circ R^\vee$ that intertwines two different 
actions on the same vector space 
\begin{equation}
  \label{Rmatr}
  R(a_1-a_2): M_1(a_1) \otimes M_2(a_2) \to M_1(a_1) \otimes_\textup{opp}
  M_2(a_2) \,. 
\end{equation}
As the word \emph{braiding} suggest, there is a constraint on the
R-matrices coming from two different ways to put three tensor factors 
in the opposite order. In our situation, the corresponding products 
of intertwiners will be simply be equal, corresponding to the
Yang-Baxter equation
\begin{multline}
R_{M_1,M_2 }(a_1-a_2) \, R_{M_1,M_3}(a_1-a_3) \, R_{M_2,M_3}(a_2-a_3) =\\=
R_{M_2,M_3}(a_2-a_3) \, R_{M_1,M_3}(a_1-a_3) \, R_{M_1,M_2}(a_1-a_2)\,, \label{YB} 
\end{multline}
satisfied by the $R$-matrices. 

\subsubsection{}\label{s_Recon}
There exists general reconstruction theorem that describe tensor categories
of a certain shape and equipped with a fiber functor to vector spaces
as representation categories of quantum groups, see \cite{Ebook1,Ebook2}. For
practical purposes, however, one may be satisfied by the following 
simple-minded approach that may be traced back to the work of 
the Leningrad school of quantum integrable systems \cite{Reshet}. 

Let $\{M_i\}$ be a collection of vector spaces over a field $\Bbbk$ and 
$$
R_{M_i,M_j}(a) \in GL(M_1\otimes M_2, \Bbbk(a)) 
$$
a collection of operators satisfying the YB equation \eqref{YB}. {}
From this data, one constructs a certain category $\cC$ of representation of
a Hopf algebra $\yY$ as follows. We first extend the R-matrices to tensor
products by the rule 
$$
R_{M_1(a_1), M_2(a_2) \otimes M_3(a_3)} = 
R_{M_1,M_3}(a_1-a_3) \, R_{M_1,M_2}(a_1-a_2) \,, \quad \textup{etc.} 
$$
It is clear that these also satisfy the YB equation. One further
extends R-matrices to dual vector spaces using inversions and
transpositions, see \cite{Reshet}. Note that in the noncocommutative 
situation, one has to distinguish between left and right dual
modules. For simplicity, we may assume that the set $\{M_i\}$ is 
already closed under duals. 

The objects of the category $\cC$ are thus $M=\bigotimes M_i(a_i)$ and 
we define the quantum group operators in $M$ as the 
as the matrix coefficients of the $R$-matrices. Concretely, we have 
\begin{equation}
\mathsf{T}_{M_0,m_0}(u) \overset{\textup{\tiny def}}= \tr_{M_0} (m_0 \otimes 1 ) R_{M_0(u), M} \in 
\End(M) \otimes \Bbbk(u) \label{trR}
\end{equation}
for any operator $m_0$ of finite rank in an auxiliary space $M_0\in
\Ob(\cC)$. The coefficients of 
$u$ in \eqref{trR} give us a supply of operators 
$$
\yY \subset \prod_{M} \End(M)\,.
$$
The YB 
equation \eqref{YB} can now be read in two different ways, 
depending on whether we designate one or two factors as 
auxiliary. With these two interpretations, it gives either: 
\begin{itemize}
\item a commutation relation between the generators \eqref{trR} of $\yY$, or 
\item a braiding of two $\yY$ modules. 
\end{itemize}
Further, any morphism in $\cC$, that is, any operator that commutes
with $\yY$ gives us relations in $\yY$ when used in the auxiliary
space. This is a generalization of the following classical fact: 
if $G \subset \prod GL(M_i)$ is a reductive algebraic group, then to
know the equations of $G$ is equivalent to knowing how $\otimes M_{k_i}$ 
decompose as $G$-modules. 

The construction explained below gives geometric $R$-matrices, that
is, geometric solutions of the YB equation acting in spaces like
\eqref{fock} in the Yangian situation, or in the 
corresponding equivariant K-theories for $\cU_\hbar(\fgh)$. This gives
a quantum group $\yY$ which is precisely of the right size for the 
enumerative application. With certain care, see \cite{MO}, the above 
construction works over a ring like $\Bbbk=\Hd_\bT(\pt)$. 

The construction of $R$-matrices uses \emph{stable envelopes}, which
is a certain technical notion that continues to find applications in 
both in enumerative and representation-theoretic contexts. 

\subsubsection{}\label{sBaxt} 
Let $z$ be an operator in each $M_i$ such that 
$$
[z \otimes z, R_{M_1,M_2}] = 0 \,. 
$$
A supply of such is provided by the Cartan torus 
$\exp(\fh)$ where, in geometric situations, $\fh$ acts by the 
ranks of the universal bundles. A classical
observation of Baxter then implies that 
$$
\left[ \mathsf{T}_{M_0,z}(u), \mathsf{T}_{M'_0,z}(u')\right] = 0 
$$
for fixed $z$ and any $M_0(u)$ and $M'_0(u')$. In general, these
depend rationally on the entries of $z$ as it is not, usually, an
operator of finite rank. The corresponding commutative subalgebras of
$\yY$ are known as the Baxter, or Bethe subalgebras. They have a 
direct geometric interpretation as the operators of quantum
multiplication, see below. 

\subsection{Stable envelopes}\label{s_Stab}

\subsubsection{}

There are two ways in which a symplectic resolution $X\xrightarrow{\,\pi\,} X_0$ may break up into 
simpler symplectic resolutions. One of them is deformation. 
There is a stratification of
the deformation space $\Def(X,\omega)=\Pic(X) \otimes_\Z \C$ by the different 
singularities that occur, the open stratum corresponding to
smooth 
$X_0$ or, equivalently, affine $X$, see \cite{Kal}. In codimension 1, one sees the simplest singularities
into which $X$ can break, and this is related to the
decompositions \eqref{pq1} and \eqref{pq2}, see \cite{BMO}, 
and so to the notion
of the K\"ahler roots of $X$ introduced above. 
In the context of Section \ref{s_mdr}, the hyperplanes of the
quantum dynamical Weyl group may be interpreted as a further 
refinement of this stratification that records
singular \emph{noncommutative} deformations of $X_0$. 

\subsubsection{}\label{s_dual}

A different way to break up $X$ into simpler pieces is to consider the 
fixed points $X^\bA$ of a symplectic torus $\bA \subset
\Aut(X,\omega)$. The first order information about the geometry of $X$ around 
$X^\bA$ is given by the $\bA$-weights in the normal bundle
$N_{X/X^\bA}$. These are called the equivariant roots for the action
of $\bA$, or just the equivariant roots of $X$ if 
$\bA$ is a maximal torus in $\Aut(X,\omega)$. 

There is a certain deep 
Langlands-like (partial) duality for equivariant symplectic
resolutions that interchanges the roles of equivariant and K\"ahler variables. The origin of this duality, 
sometimes called the 3-dimensional mirror symmetry, is in 3-dimensional 
supersymmetric gauge theories, see \cite{IS,BHOO} and also e.g.\ \cite{BDG,BDGH, BLPW1} for
a thin sample of references. As this duality interchanges K-theoretic 
enumerative information, the quantum difference equation for $X$ 
becomes the shift operators for its dual $X^\vee$, see \cite{PCMI} for an 
introduction to these notions. Among other things, 
the K\"ahler and equivariant roots 
control the poles in these difference equations, which makes it 
 clear that they should be exchanges by the duality. 

For the development of the theory, it is very important to be able to
both break $X$ into simpler piece and also to find $X$ as such a piece in
a more complex ambient geometry. 

\subsubsection{} 

For instance for $X=\Hilb(\C^2,n)$, we can take the maximal 
torus $\bA \subset
SL(2)$, in which case $X^\bA$ is a finite set of monomial ideals 
$$
I_\lambda \subset \C[x_1,x_2]\,, \quad |\lambda|=n \,, 
$$
indexed by partitions $\lambda$ of the number $n$.  At $I_\lambda$, the
normal weights 
are $\{\pm \textup{hook}(\square)\}_{\square\in \lambda}$,
with the result that
\begin{equation}
\textup{equivariant roots of  } \Hilb(\C^2,n)  = \{\pm 1,\dots,\pm n\} 
\,. \label{ER_hilb}
\end{equation}
Notice the parallel with \eqref{KR_hilb}. 
It just happens that $\Hilb(\C^2,n)$ is 
self-dual, we will see other manifestations of this below. 

\subsubsection{}\label{sMrn}

More importantly, products $\prod_{i=1}^r \Hilb(\C^2,n_i)$ may be
realized in a very nontrivial way as fixed loci of a certain torus
$\bA$ on an ambient variety $\cM(r,\sum n_i)$. Here 
$\cM(r,n)$ is the moduli space of framed torsion free sheaves $\cF$ on 
$\bP^2$ of rank $r$ and $c_2(\cF)=n$. A framing is a choice of an 
isomorphism 
$$
\phi: \cF\big|_{L}  \xrightarrow{\,\, \sim \,\,} 
\cO_{L}^{\oplus r} \,, \quad L = \bP^2 \setminus \C^2 \,, 
$$
on which the automorphism group $GL(r)$ acts by postcomposition. 
The spaces $\cM(r,n)$ are the general Nakajima varieties associated to 
the quiver with one vertex and one loop. They play a
central role in supersymmetric gauge theories as symplectic
resolutions of the Uhlenbeck spaces of framed instantons, see
\cite{NakL}. It is easy to see that $\cM(1,n)= \Hilb(\C^2,n)$ and 
$$
\cM(r,n)^{\bA} = \bigsqcup_{\sum n_i=n} \prod_{i=1}^r \Hilb(\C^2,n_i)
\,, 
$$
where 
$\bA 
\subset GL(r)$ is the maximal torus. For general Nakajima 
varieties, there is a similar decomposition for the maximal torus 
$\bA$ of framing automorphisms. 

\subsubsection{}
The cohomological stable envelope is a certain 
Lagrangian correspondence 
\begin{equation}
\Stab \subset X \times X^\bA \,,\label{Stab1}
\end{equation}
which may be seen as an improved version of the attracting manifold 
$$
\Attr = \left\{(x,y) \big| \lim_{a\to 0} a\cdot x =y\right\} \subset
X\times X^\bA \,. 
$$
The support of \eqref{Stab1} is the full attracting set $\Attr_f \subset
  X\times X^\bA$, which is the smallest closed subset 
that contains the diagonal and is closed under taking 
$\Attr(\,\cdot\,)$.

To define attracting and repelling
manifolds, we need to separate the roots for $\bA$-weights into positive 
and negative, that is, we need to choose a chamber  $\fC\subset \Lie \bA$ 
in the complement of the root hyperplanes.  Note this gives an ordering on the set of 
components $\bigsqcup F_i = X^\bA$ of the fixed locus: $F_1 > F_2$ if 
$\Attr_f(F_1)$ meets $F_2$.

For 
\begin{equation}
\Lie \bA = \{ \diag(a_1,\dots,a_r) \} \subset \gl(r) \label{LieA}
\end{equation}
as in Section \ref{sMrn}, the roots are $\{a_i-a_j\}$, and so a choice
of $\fC$ is the usual choice of a
Weyl 
chamber. As we will see, it will 
correspond to an ordering of tensor factors as in Section
\ref{s_Rmat}.

The stable envelope is
characterized by:
\begin{itemize}
\item it is supported on $\Attr_f$, 
\item it equals\footnote{The $\pm$ a choice here is reflecting a choice
    of a \emph{polarization} of $X$, which is a certain auxiliary
    piece of data that one needs to fix in the full development of the theory.} $\pm \Attr$ near the diagonal in 
$X^\bA \times X^\bA  \subset X \times X^\bA$, 
\item for an off-diagonal component $F_2 \times F_1$ of $X^\bA \times
  X^\bA $, we have
  \begin{equation}
\deg_{\Lie(\bA)} \Stab\Big|_{F_2 \times F_1} < \tfrac12 \codim F_2 =
\deg_{\Lie(\bA)} \Attr\Big|_{F_2 \times F_2}\,, \label{degStab}
\end{equation}
where the degree is the usual degree of polynomials for 
$\Hd_\bA(X^\bA,\Z) \cong \Hd(X^\bA) \otimes \Z[\Lie \bA]$. 
\end{itemize}

Condition \eqref{degStab} is 
a way to quantify the idea that $\Stab\big|_{F_2\times F_1}$ is smaller than $\Attr(F_2)$.  This makes the 
stable envelope a canonical representative of
$\overline{\Attr(F_1)}$ modulo cycles supported on the lower 
strata of $\Attr_f$.

The existence and uniqueness of stable envelopes are proven, under very general
assumptions on $X$ in \cite{MO}. As these correspondences are canonical,
they are invariant under the centralizer of $\bA$ and, in particular,
act in $\bT$-equivariant cohomology for any ambient torus $\bT$. 

\subsubsection{}

To put ourselves in the situation of Section \ref{s_Recon}, we define a
category in which the objects are 
\begin{equation}
\bF(a_1,\dots,a_r) = \bigoplus_{n\ge 0} \Hd_\bT(\cM(r,n))\label{bF}
\end{equation}
and the maps defined by stable envelopes, like 
\begin{equation}
  \label{Stabpm}
  \xymatrix{\bF(a_1) \otimes \bF(a_2) \ar@/^1pc/[rrr]^{\Stab_+} 
\ar@/_1pc/[rrr]_{\Stab_-}  &&&  \bF(a_1,a_2) }\,,
\end{equation}
where $\fC_{\pm} = \{a_1 \gtrless a_2\}$, are declared to be
morphisms. Since both maps in \eqref{Stabpm} are isomorphisms 
after $\bA$-equivariant localization, we get a rational matrix
\begin{equation}
R(a_1-a_2) = \Stab_-^{-1} \circ \Stab_+  \in \End(\bF(a_1) \otimes
\bF(a_2)) \otimes \Q(a_1-a_2) \,.\label{RM}
\end{equation}
Its basic properties are summarized in the following 

\begin{Theorem}[\cite{MO}] The $R$-matrix \eqref{RM} 
satisfies the YB equation
  and defines, as in Section \ref{s_Recon}, an action of $\yY(\glh(1))$ in
  equivariant cohomology of $\cM(r,n)$. The Baxter subalgebras in
  $\yY(\glh(1))$ are the algebras of operators of quantum
  multiplication. In particular, the vacuum-vacuum elements of the
  $R$-matrix are the operators of classical multiplication in
  $\cM(r,n)$. 
\end{Theorem}

Here $z$ in the Cartan torus of $\glh(1)$ acts by $z^n$ in the $n$th
term of \eqref{bF}, which clearly commutes with $R$-matrices. The 
vacuum in \eqref{bF} is the $n=0$ term. 

\subsubsection{}
For a general Nakajima variety, it is proven in \cite{MO} that the 
corresponding R-matrices define an action of $\yY(\fg)$ for a 
certain Borcherds-Kac-Moody Lie algebra 
$$
\fg = \fh \oplus \bigoplus_{\alpha\ne 0} \fg_\alpha \,,
$$
with finite-dimensional root spaces $\fg_\alpha$. This Lie algebra is 
additionally graded by the cohomological degree and it has been 
conjectured in \cite{OkKac} that graded dimensions of $\fg_\alpha$ are 
given by the Kac polynomial for the dimension vector $\alpha$. 
A slightly weaker version of this conjecture is proven in \cite{SV2,SV3}. 

Again, the operators of classical multiplication are given by the 
vacuum-vacuum matrix elements of the R-matrix, while for the quantum
multiplication we have the formula already announced in Section \ref{s_pq}

\begin{Theorem}[\cite{MO}] For a general Nakajima variety, 
quantum multiplication by divisors is given by the formula \eqref{pq2}
and hence the quantum differential equation is the Casimir connection 
for $\yY(\fg)$. 
\end{Theorem}

\subsubsection{}
Back to the Hilbert scheme case, the $R$-matrix \eqref{RM} acting in
the tensor product of two Fock spaces is a very important object for
which various formulas and descriptions are available. 
The following description was obtained in \cite{MO}. 

The operators 
$$
\alpha^{\pm}_n = \alpha_n \otimes 1 \pm 1 \otimes \alpha_n 
$$
act in the tensor product of two Fock spaces, and form two commuting 
Heisenberg subalgebras. They are analogous to the center of mass and
separation coordinates in a system of two bosons, and we 
can similarly decompose 
\begin{equation}
\bF(a_1) \otimes \bF(a_2) = \bF_{+}(a_1+a_2) \otimes
\bF_{-}(a_1-a_2)\,,\label{Fpm}
\end{equation}
where to justify the labels we introduce the zero modes 
$\alpha_0(\gamma)$ that act on $\bF(a_i)$ by $-a_i \int \gamma$. 
Here integral denotes the equivariant integration of $\gamma \in 
\Hd_\bT(\C^2)$. 
Consider the operators $L_n$ defined by 
\begin{equation}
\sum L_n \zeta^{-n} = \frac14 \, \textbf{:} \bal_-(\zeta)^2 \textbf{:}  \, 
\pm \frac12 \, \hbar \, \partial \bal_-(\zeta) - \frac14 \int \hbar^2
\,, 
\label{FFVir}
\end{equation}
where $\textbf{:} \bal_-(\zeta)^2 \textbf{:}$ is the normally ordered 
square of the operator $\bal_-(\zeta)= \sum \alpha^{-}_n \zeta^n$,
which now has a constant term in $\zeta$, and 
$\partial$ stands for $\zeta \frac{\partial}{\partial \zeta}$. For the 
cohomology arguments of $\alpha^-_n$, we use the conventions of Section
\ref{s_Lehn}. For either choice of sign, \eqref{FFVir} form 
the Virasoro algebra 
\begin{equation}
  \label{commL}
  [L_n,L_m] = (m-n) L_{n+m} + \frac{(n^3-n)}{12} \left(
1+ 6 \int \hbar^2 \right) \delta_{n+m} \,, 
\end{equation}
in a particular free field realization that is very familiar from the 
work of Feigin and Fuchs \cite{FeiginFuchs} 
and from CFT. Note 
that, for generic $a_1-a_2$, $\bF_-(a_1-a_2)$ is irreducible with highest weight 
$$
L_0 \, \big|\big\rangle = \frac14 \, \int ((a_1-a_2)^2 - \hbar^2) \, 
\big|\big\rangle \,, 
$$
which, like \eqref{commL} is invariant under $a_1 \leftrightarrow 
a_2$ and $\hbar \mapsto - \hbar$

\begin{Theorem}[\cite{MO}]\label{tR} 
The $R$-matrix is the unique operator in 
$\bF_-$ in \eqref{Fpm} that 
preserves the vacuum $\big|\big\rangle$ and 
interchanges the two signs in \eqref{FFVir}. 
\end{Theorem}

This is the technical 
basis of numerous fruitful connection between $\yY(\glh(1))$
and CFT. For instance, the operator $R^\vee = (12) \circ R$ is related to
 the reflection operator in Liouville CFT, see \cite{ZZ},
and the YB equation satisfied by $R$ reveals new unexpected
features of affine Toda field theories. 

\section{Some further directions}

\subsection{K-theoretic counts}\label{s_Kcounts}

\subsubsection{} 
K-theoretic counts require a definite technical investment to
be done properly, but offer an ample return producing 
deeper and more symmetric theories. 
For example, the dualities 
already mentioned in Section \ref{s_dual} reveal their full power only in 
equivariant K-theory. 

Focusing on $\Hilb(\C^2,n)$, its self-duality may be 
put into an even deeper and, not surprisingly, widely conjectural framework of M-theory,
which is a certain unique 11-dimensional supergravity theory with the 
power to unify many plots in modern theoretical physics as well as 
pure mathematics.  Its basic actors are membranes with a 3-dimensional 
worldvolume that may appear as strings or even point particles to a
low-resolution observer, see \cite{Baggeretal} for a recent review. The
supersymmetric index for membranes of the form
$$
C \times S^1 \subset  Z \times S^1\,,
$$
where $C$ is a complex curve in a 
complex Calabi-Yau 5-fold $Z$, should be a virtual representation of
the 
automorphisms $\Aut(Z,\Omega^5_Z)$ that preserve the 5-form, so a certain K-theoretic curve count in $Z$. 
See \cite{MDT} for what it might look like and a conjectural 
equivalence with K-theoretic DT counts for 
$$
Y = Z^{\Ct_z}\,, 
$$
for any $\Ct_z \subset \Aut(Z,\Omega^5_Z)$ with a purely 3-dimensional
fixed locus. The most striking feature of this equivalence is the 
interpretation of the variable $z$ 
\begin{align}
z & = \textup{degree-counting variable in \eqref{starLehn}} 
\notag \\
& = \textup{the $\chi(\cO_C)$-counting variable in DT theory, see 
Section \ref{szn}}  \notag  \\
&\overset{!}= \textup{equivariant variable $z\in \Ct_z$ in M-theory}
  \,. \label{Mz} 
  \end{align}
Note that in cohomology equivariant variables take values in a Lie
algebra, and so cannot be literally on the same footing as K\"ahler
variables. 

\subsubsection{}

As a local model, one can take the total space 
\begin{equation}
Z = 
\begin{matrix}
\cL_1 \oplus \cL_2 \oplus \cL_3 \oplus \cL_4 \\
\downarrow\\
B
\end{matrix} \,,\quad \bigotimes_{i=1}^4 \cL_i = \cK_B\,, 
\label{local_curveM}
\end{equation}
of 4 line bundles over a smooth curve $B$. In this geometry, one can 
designate any two $\cL_i$ to form $Y$ by making $z$ scale the other 
two line bundles with opposite weights. Thus the counting of Section 
\ref{sTQFT} with $\sS=\C^2$, when properly set up in K-theory,
 has a conjectural $S(4)$-symmetry that permutes the weights 
$$
(t_1,t_2,\frac{z}{\sqrt{q t_1 t_2 }}, \frac{1}{z \sqrt{q t_1 t_2}}) \,,
\quad 
\begin{pmatrix}
  t_1 \\ & t_2 
\end{pmatrix} \in \Aut(\C^2)  \,, 
$$
where $q^{-1}$ is the Chern root of $\cK_B$.  

\subsubsection{}

This theory is described by certain $q$-difference equation that
describes the change of the counts as the degrees of $\{\cL_i\}$
change, and correspondingly the variables $t_1,t_2,z$ are shifted 
by powers of $q$. In fact, similar $q$-difference equations can be 
defined for any Nakajima variety including $\cM(r,n)$ from Section
\ref{sMrn}. For those, there are also difference equations in the 
framing equivariant variables $a_i$. 

\begin{Theorem}[\cite{PCMI,OS}]\label{t_qKZ} The $q$-difference equations in 
variables $a_i$ are the quantum Knizhnik-Zamolodchikov (qKZ) equations 
for $\cU_\hbar(\glhh(1))$ and the $q$-difference equation in $z$ is 
the lattice part of the dynamical Weyl group of this 
quantum group. Same is true for a general Nakajima variety and 
the corresponding quantum loop algebra $\cU_\hbar(\fg)$. 
\end{Theorem}

\noindent 
The qKZ equations, introduced in \cite{FrenResh}, 
have the form 
\begin{equation}
  \label{qKZ} 
  \Psi(q a_1, a_2,\dots, a_n) = ( z \otimes 1 \otimes \dots \otimes 1)
  \, 
R_{1,n}(a_1/a_n) \dots R_{1,2}(a_1/a_2) \, \Psi 
\end{equation}
with similar equations in other variables $a_i$, where $z$ is as in
\eqref{LieA}. For the $R$-matrices of $\cU_\hbar(\fgh)$, where 
$\dim \fg < \infty$, these play the same role
in integrable 2-dimensional lattice models as the classical KZ 
equations play in their conformal limits, see e.g.\ \cite{JMbook}. 

For 
$\cU_\hbar(\glhh(1))$, the $R$-matrix is a generalization of the 
$R$-matrix of Theorem \ref{tR} constructed using the K-theoretic 
stable envelopes. The algebra $\cU_\hbar(\glhh(1))$ is constructed
from this R-matrix using the general procedure of Section \ref{s_Recon}. 
For this particular geometry, it also coincides with the algebras 
constructed by many authors by different means, including explicit 
presentations see e.g.\  \cites{BurSch, FeiJ1,Negut, SV1} for a sample of references where the 
same algebra appears. 

It takes a certain development of the theory to define and make 
concrete the operators from the quantum dynamical Weyl group 
of an algebra like $\cU_\hbar(\fgh)$. They crucially depend on 
certain features of K-theoretic stable envelopes that are not 
visible in cohomology. 

\subsubsection{}\label{s_StabK}
In the search for the right K-theoretic generalization of
\eqref{degStab}, one should keep in mind that the correct notion of the degree
of a multivariate Laurent polynomial is its Newton polytope, that is,
the convex hull of its exponents, considered \emph{up to
  translation}. Translations correspond to multiplication by
monomials, which are invertible functions. The right generalization 
of \eqref{degStab} is then 
  \begin{equation}
\textup{Newton}\left(\Stab\Big|_{F_2 \times F_1}\right)  \subset 
\textup{Newton}\left(\Attr\Big|_{F_2 \times F_2}\right) + 
\textup{shift}_{F_1,F_2}\,, \label{degStabK}
\end{equation}
for a certain collection of shifts in $\bA^\vee \otimes \R$, where 
$\bA^\vee$ is the character lattice of the torus $\bA$. Shifts come 
from weights of fractional line bundles $\cL \in \Pic(X) \otimes \R$ 
at the fixed points, see \cite{PCMI} for a survey. This fractional line 
bundle, called the \emph{slope} of the stable envelope is a new 
parameter, the dependence on which is locally constant and quasiperiodic
--- if the slopes differ by integral line bundle $\cL$, then the
corresponding stable envelopes differ by a twist by $\cL$. 

Stable envelopes change across certain hyperplanes in 
$\Pic(X) \otimes \R$ that form a $\Pic(X)$-periodic hyperplane 
arrangement closely related to the K\"ahler roots of $X$.  For
example, for $\glh(1)$ the roots are $\{\alpha\} = \N\setminus \{0\}$ and so the 
affine root hyperplanes are all rationals 
\begin{equation}
  \{x  | \, \exists \alpha  \, \langle \alpha, x \rangle \in \Z\}  = 
\Q \subset \R\,. 
\label{wallsQ}
\end{equation}
The dynamical 
Weyl group element $\bB_{a/b}$ in \eqref{Bwall} 
corresponding to a wall 
$\frac ab \in \Q$ acts trivially in $K_\bT(\Hilb(\C^2,n)$ if $b > n$, 
reflecting the fact that the K\"ahler roots of $\Hilb(\C^2,n)$
form a finite subset 
$\{\pm 1, \dots, \pm n \}$ of the roots of 
$\glh(1)$.

The change 
across a particular wall is recorded by a certain wall $R$-matrix 
$R_\textup{wall}$ and the $R$-matrices corresponding to a change of 
chamber $\fC$ as in \eqref{RM}
 factor into an infinite product of those, see \cite{PCMI,OS}. 
Each term in such factorization corresponds to a certain root subalgebra 
$$
\cU_\hbar(\fg_\textup{wall}) \subset \cU_\hbar(\fgh)
$$
stable under the action of 
\begin{equation}
\widehat{z} \in \textup{Cartan torus of $\fg$} \times \Ct_q
\,,\label{hatz}
\end{equation}
where $\Ct_q$ acts by loop rotation automorphisms of
$\cU_\hbar(\fgh)$. One defines \cite{OS} the dynamical Weyl group by taking certain 
specific $(z,q)$-dependent elements
\begin{equation}
\bB_\textup{wall} \in \cU_\hbar(\fg_\textup{wall})  \,. \label{Bwall}
\end{equation}  
For finite-dimensional $\fg$, we have $\fg_\textup{wall} =
\mathfrak{sl}(2)$ for every wall and the construction specializes to
the classical construction of Etingof and Varchenko \cite{EV}. The
operators $\bB_\textup{wall}$ satisfy the braid relations of the wall 
arrangement. Because each $\bB_\textup{wall}$ depends on 
\eqref{hatz} through the equation of the corresponding wall, these relations look like the Yang-Baxter equations \eqref{YB}, in which each term 
depends on $a=(a_1,a_2,a_3)$ through the equation $a_i-a_j=0$ of the 
hyperplane being crossed in $\Lie \bA$. 

\subsection{Monodromy and derived equivalences}\label{s_mdr}

\subsubsection{}
As a special $\widehat{z}$-independent case, the dynamical Weyl group
contains the so-called: 
\begin{itemize}
\item[\bo{w}]  quantum Weyl group of $\cU_\hbar(\fgh)$, which
plays many other roles, including 
\item[\bo{m}] this is the monodromy group of the quantum differential
  equation
\eqref{qde}, and  
\item[\bo{p}] this group describes the action on $K_\bT(X)$ of the derived
  automorphisms of $X$ constructed by Bezrukavnikov and Kaledin using 
quantizations $\hX_c$ of $X$ 
in characteristic $p\gg 0$, see \cites{BezF,BK2,BezL,Kaledin_q,Los1,Los2,Los3,Los4}. 
It thus plays the same role in modular representation theory of 
$\hX_c$ as the 
Hecke algebra plays in the classical Kazhdan-Lusztig theory. 
\end{itemize}

\begin{Theorem}[\cite{BO}] We have 
$\textup{\bo{w}}=\textup{\bo{m}}=\textup{\bo{p}}$ for all Nakajima varieties. 
\end{Theorem}

Other known infinite series of equivariant symplectic
resolutions are also considered in \cite{BO}. For finite-dimensional $\fg$, 
the description of the monodromy 
of the Casimir connection via the quantum Weyl group is a conjecture
of V.~Toledano-Laredo 
\cite{TolCas}. The equality $\textup{\bo{m}}=\textup{\bo{p}}$
for all equivariant symplectic resolutions is a conjecture of
Bezrukavnikov and the author, see the discussion in
 \cite{ABM,slc}. 

Here $\hX_c$ is an associative algebra deformation of the Poisson
algebra of functions on $X$ or $X_0$. While it may be studied
abstractly, quantizations of Nakajima varieties 
may be described concretely as 
quantum Hamiltonian reductions, see \cite{EtCal,GanGinz}. For 
$X=\Hilb(\C^2,n)$, $\hX_c$ is the algebra generated by 
symmetric polynomials in $\{w_1,\dots,w_n\}$, and the operators 
of the rational Calogero quantum integrable systems --- a commutative
algebra of differential operators that includes the following rational 
analog of \eqref{HCS}
\begin{equation}
H_{\textup{C,rat}} = \tfrac12 \sum_{i=1}^N \frac{\partial^2}{\partial
    w_i^2}  - c(c+1) \sum_{i<j\le N}
\frac{1}{(w_i-w_j)^2} \,. \label{HCrat}
\end{equation}
The kinship between \eqref{HCS} and \eqref{HCrat} is closer than
normal because of the self-duality of $\Hilb(\C^2,n)$. Recall that 
duality swaps equivariant and K\"ahler variables and instead of an 
equivariant variable $\theta$ in \eqref{HCS} we have
a K\"ahler variable $c$  in \eqref{HCrat}. It parametrizes deformations of $\hX_c$ in 
the same way as $\Pic(X)\otimes_\Z \Bbbk$ parametrizes 
deformations of $X$ over a field $\Bbbk$. 

For $p\gg 0$, the BK theory produces derived equivalences 
\begin{equation}
  \label{BK}
\Db \Coh X^{(1)} \xleftrightarrow{\quad\sim\quad}
\Db \hX_c\textup{-mod} 
\xleftrightarrow{\quad\sim\quad}
\Db \Coh X^{(1)}_\textup{flop}
\end{equation}
for every nonsingular value of $c\in\Z$. A shift $c\mapsto c+p$ 
twists \eqref{BK} by $\cO(1)$. Here $X^{(1)}$ denotes the 
Frobenius twist of $X$ and $X_\textup{flop}$ refers to a change of 
stability condition in the GIT construction of $X$. 

In the  $\textup{\bo{m}}=\textup{\bo{p}}$ interpretation, 
the composed equivalence in \eqref{BK} becomes the transport
of the QDE from the point $z=0$, that corresponds to $X$, to the 
point $z=\infty$, that corresponds to $X_\textup{flop}$, along the ray
with 
$$
\arg z = - 2 \pi \, \frac{c}{p} \,. 
$$
In particular, this identifies the singularities of the QDE, given by
the K\"ahler roots of $X$, with the limit $\lim_{p\to\infty} \frac1p
\{c_\textup{sing}\}$ of the singular parameters of the quantization. 
See \cite{BO} for details. 

\subsubsection{}

Monodromy of a flat connection with regular singularities is an analytic map between algebraic 
varieties that may be seen as a generalization of the exponential map
of a Lie group. There is a long tradition, going back to at least the 
work of Kohno and Drinfeld of computing the monodromy of 
connections of representation-theoretic origin in terms of closely  
related algebraic structures. Just like for the exponential map, there
is a certain progression in this, as one goes from additive variables 
to multiplicative, and also from multiplicative --- to
elliptic. E.g.\ in the case at hand the QDE (= the Casimir
connection) is defined for modules $M$ over the Yangian $\yY(\fg)$, and is 
computed in terms of the action of $\cU_\hbar(\fgh)$ in a closely
related representation\footnote{geometrically, the right relation
  between $H_\bT(X)$, where the Yangian acts, and $K_\bT(X)$ of 
$X$, where $\cU_\hbar(\fgh)$ acts,  is given 
by a certain $\Gamma$-function analog of the Mukai vector, as in the 
work of Iritani \cite{Ir}}, see in particular \cite{GauTL}. 

A key step in capturing the monodromy algebraically is typically a 
certain 
compatibility constraint between the monodromy for $M=M_1\otimes
M_2$ and the monodromy for the tensor factors. The framework
introduced above gives a very conceptual and powerful way to prove such
statements. Recall that, geometrically, $\otimes$ 
arises as a special correspondence between $X^\bA$ and $X$, where 
$\bA$ is a torus that acts on $X$ preserving the symplectic form. 
Therefore, it is natural to ask, more generally, for a compatibility
between the monodromy of the QDE for $X$ and $X^\bA$. 

In fact, one can ask a more general question about the compatibility
of the corresponding $q$-difference equations as in Theorem \ref{t_qKZ}.  
Let 
$$
\bZ = \Pic(X) \otimes_\Z \Ct
$$
be the Kahler torus of $X$ and $\bZb$ be the toric compactification 
of $Z$ corresponding to the fan of ample cones of flops of $X$. 
Its torus-fixed points $0_{X_\textup{flop}}\in \bZb$ correspond to all possible flops
of $X$. 
A regular $q$-difference connection on a smooth toric 
variety  is an action of the the cocharacter lattice 
$$
\Pic(X) \owns \lambda \to q^{\lambda} \in \bZ \,, \quad q\in \Ct \,, 
$$ 
on a 
(trivial) vector bundle over $\bZb$. Shift operators define a
commuting regular $q$-difference connections in the variables 
$a \in \bA \subset \bAb$, where $\bAb$ is the toric variety given by 
the fan of the chambers $\fC$. The $q$-difference connection for $X^\bA$ 
sits overs the torus fixed points $0_\fC\in \bAb$. 

The most interesting analytic feature here is that the connections in
$z$ and $a$, while compatible and separately regular, are \emph{not} regular
jointly. This can never happen for differential equations, see
\cite{Del2}, but is commonplace for $q$-difference equations 
as illustrated by the system: 
$$
f(qz,a) = a f(z,a) \,, \quad f(z,qa) = z f(z,a) \,. 
$$
As a result, near any point $(0_X,0_\fC) \in \bZb \times \bAb$, we get 
two kinds of solutions. Those naturally arising enumeratively are holomorphic in $z$ in a punctured 
neighborhood of $0_X$ and meromorphic in $a$ with poles accumulating 
to $0_\fC$. These may be called the $z$-solutions. For $a$-solutions, 
the roles of $z$ and $a$ are exchanged. These naturally appear in the 
Langlands dual setup and the initial conditions at $a=0_\fC$ for them
are the $z$-solutions for $X^\bA$. 

Transition matrices between the $a$-solutions and the
$z$-solutions, which is by construction elliptic, intertwine
the monodromy for $X^\bA$ and $X$, and vice versa. Note 
these transition matrices may, in principle, be computed from the
series expansions near $(0_X,0_\fC)$, which differentiates them from
more analytic objects like monodromy or Stokes matrices. 

\begin{figure}[!htbp]
  \centering
   \includegraphics[scale=0.75]{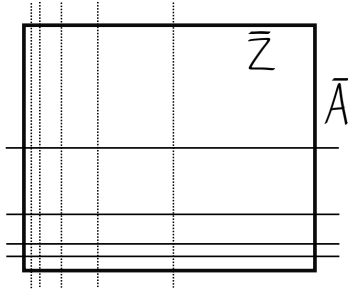}
 \caption{$z$-solutions are convergent power series in $z$ with
   coefficients in $\Q(a)$, and the poles (solid lines in the picture)
   of these coefficients
   accumulate to $a=0_\fC$. The poles of $a$-solutions, the dashed
   lines in the figure, similarly accumulate to $z=0_X$.}
\label{f_az}
\end{figure}

\begin{Theorem}[\cite{AO}]
The transformation from the $a$-solutions to $z$-solutions is given by 
elliptic stable envelopes, a certain elliptic analog of $\Stab$
constructed in \cite{AO}. 
\end{Theorem}

\subsubsection{}

To complete the picture, one can give Mellin-Barnes-type 
integral solutions to the quantum $q$-difference equations, 
and so, in particular to the qKZ and dynamical equations 
for $\cU_\hbar(\fgh)$ in tensor products of evaluation 
representations \cite{AObethe}.  It is well-known that the stationary phase 
$q\to 1$ limit in such integrals diagonalizes the Baxter 
subalgebra and hence generalizes the classical ideas of 
Bethe Ansatz.

\vfill 

\noindent 
Andrei Okounkov\\
Department of Mathematics, Columbia University\\
New York, NY 10027, U.S.A.\\

\vspace{-12 pt}

\noindent 
Institute for Problems of Information Transmission\\
Bolshoy Karetny 19, Moscow 127994, Russia\\

\vspace{-12 pt}

\noindent 
Laboratory of Representation
Theory and Mathematical Physics \\
Higher School of Economics \\ 
Myasnitskaya 20, Moscow 101000, Russia 

\end{document}